\documentstyle [12pt,leqno]{article}

\def\P#1{\mbox{${\bf P}^{#1}$}}
\def\CM{Cohen-Macaulay}
\def\aCM{arithmetically Cohen-Macaulay}
\def\O{\mbox{${\cal O}$}}
\def\QED{\begin{flushright} QED \end{flushright}}
\def\M{{\sc Macaulay}}

\newtheorem{fact}{Proposition}

\begin{document}

\quad
\bigskip

\begin{center}
{\bf Arithmetically Cohen-Macaulay Curves cut out by  Quadrics}

Sheldon Katz

Department of Mathematics

Oklahoma State University

Stillwater, OK 74078
\end{center}

The following question was raised by M. Stillman.

\bigskip\noindent
{\bf Main Question:} Let $C\subset\P r={\bf C}\P r$ be a smooth arithmetically
Cohen-Macaulay
curve which is cut out scheme theoretically by quadrics.  Is the homogeneous
ideal of $C$ necessarily cut out by quadrics?

\bigskip
In \cite{EEK}, it was shown that the question has an affirmative answer if
$r\le
5$.  The purpose of this note is to show that the question has a
negative answer (there is a counterexample with $r=7$).

It is a pleasure to thank D. Eisenbud, M. Stillman, and C. Walter for helpful
conversations.

\section{Homogeneous and Scheme-Theoretic Generation by Quadrics}

Let $X$ be a projective variety.  It is often of interest to know whether or
not
the homogeneous ideal of $X$ can be generated by quadrics, e.\hskip0em g.\ if
$X$ is a
general canonical curve.  In such a case, $X$ is cut out scheme-theoretically
by
quadrics as well.  It is usually easier to verify the scheme-theoretic
statement--- this amounts to ignoring the vertex of the affine cone over $X$.

\bigskip\noindent
{\bf Problem:} Let $C\subset\P r={\bf C}\P r$ be a smooth curve which is cut
out scheme theoretically by quadrics.  Is the homogeneous
ideal of $C$ necessarily cut out by quadrics?

\bigskip

In \cite{EEK}, this problem was investigated.  The answer is a resounding
\underbar{no}.  A counterexample was found with $r=5$.  However, positive
results
were found.  The problem has an
affirmative answer for curves on scrolls, all curves with $r\le 4$, and
arithmetically Cohen-Macaulay curves which lie on projectively normal {\rm K}3
surfaces cut out by quadrics (this includes all arithmetically C-M curves with
$r=5$).  This leads to a more precise question, which we could not answer:

\bigskip\noindent
{\bf Question:} Let $C\subset\P r={\bf C}\P r$ be a smooth arithmetically
Cohen-Macaulay
curve which is cut out scheme theoretically by quadrics.  Is the homogeneous
ideal of $C$ necessarily cut out by quadrics?

\bigskip

It turns out that this question also has a negative answer.

\begin{fact}
\label{main}
Let $C\subset\P 7$ be a general degree 19 embedding of a general genus 12 curve
over an algebraically closed field of characteristic 0.  Then $C$ is smooth and
arithmetically Cohen-Macaulay, $C$ is cut out scheme-theoretically by quadrics,
and the homogeneous ideal of $C$ is not cut out by quadrics.
\end{fact}

\section{Candidates for a Counterexample}

Let $C\subset\P r$ be an arithmetically \CM\ curve of degree $d$ and genus $g$.
Assume in addition that $\O_C(1)$ is non-special, i.\hskip0em e.\
$H^1(\O_C(1))=0$.  Then $d=g+r$.

It turns out that for certain values of $g$ and $r$, the homogeneous ideal of
such a curve $C$ {\em cannot\/} be cut out by quadrics, for simple dimension
reasons.  Let $I$ denote the ideal sheaf of $C$.  Then if

\begin{equation}
\label{ineq23}
(r+1)h^0(I(2))<h^0(I(3)),
\end{equation}

the natural map

\begin{displaymath}
H^0(I(2))\otimes H^0(\O_{\P r}(1))\to H^0(I(3))
\end{displaymath}
{\em cannot\/} be surjective, so that the homogeneous ideal of $C$ cannot be
generated by quadrics.  Using Riemann-Roch, (\ref{ineq23}) becomes
\begin{equation}
\label{glow}
 g>\frac{r(r-2)}3.
\end{equation}

On the other hand, if we want $C$ to be scheme-theoretically cut out by
quadrics,
then we must have enough quadrics, i.\hskip0em e.

$$\pmatrix{r\cr
        2}-g\ge r-1.$$

Equality holds if and only if $C$ is a complete intersection of $r-1$ quadrics;
but in this case the homogeneous ideal is cut out by quadrics as well.
This can be improved slightly: in \cite[Cor.\ 2.5]{EEK} it was shown that if
$C$
is cut out scheme theoretically by $r$ quadrics, then necessarily

\begin{equation}
\label{r2}
g=(r-1)d/2+
1-2^{r-1}.
\end{equation}
So if (\ref{r2})
does not hold, then

\begin{equation}
\label{gup}
g\le\frac{r^2-3r-2}2
\end{equation}

There are no counterexamples to the main question for $r\le 5$ \cite{EEK}.
Suppose that there is a non-special
counterexample with $r=6$. Then $g\ge 9$
by (\ref{glow}).  Since (\ref{gup}) gives $g\le 8$, it follows that (\ref{r2})
holds,
and $g=9$.  But then $d=15$, and a contradiction is reached.

Turning next to $r=7$, (\ref{glow}) gives $g\ge 12$, and (\ref{gup}) gives
$g\le 13$.
In the following section, we show that in fact the {\em general \/} curve of
degree 19 and genus 12 in \P7 is a counterexample.

\section{The counterexample}

Pick 22 general points $p_1,p_2,p_3,q_1,\ldots ,q_7,r_1,\ldots ,r_{12}$ in
$\P2$.
Let $C'$ be a general plane curve of degree 9 passing through the $p_i$ with
multiplicity 3, through the $q_i$ with multiplicity 2, and simply through the
$r_i$.  The linear system $|L|$ of degree 7 curves passing doubly through the
$p_i$ and
simply through the $q_i$ and $r_i$ maps $C'$ birationally to a smooth
curve $C$ of
degree 19 and arithmetic genus 12 in $P^7$.

It is a simple matter to use \M\ \cite{BS} to construct such a curve.  In
describing the calculation, I will informally say that a general curve has a
certain property, when I mean that the property is satisfied for an example
curve constructed using \M's pseudo-random number generator.  In fact, I
repeated
the construction several times with different pseudo-random coefficients, and
the properties mentioned below held in each instance.  Thus, as expected, a
``general'' curve has been constructed.

\M's pseudo-random number generator is used to construct 22 ``general'' points
in
$\P2_{{\bf F}_{31991}}$, and from this the
curve $C'$ (actually, there is no harm in supposing that the $p_i$ are
$(1,0,0),
(0,1,0),(0,0,1)$, to shorten computations).  By calculating the Jacobian of
$C'$, it is checked that the singular
scheme of $C'$ has degree 19 as expected  (triple points count at least 4
times).
Hence
$C'$ has the expected geometric genus 12.  The equations of the image curve
$C$ can then be explicitly
calculated.  $C$ is cut out ideal theoretically by 9 independent quadrics and 2
independent cubics, and has Hilbert function $(1+6t+12t^2)(1-t)^{-2}$.  In
particular $C$ has arithmetic genus 12; being the image of the normalization of
$C'$ by the base point free system $|L|$ on the blowup of \P2, it follows that
$C$ is smooth.  Let
${\tilde C}\subset \P2$ be the scheme cut out by the 9
quadrics
alone.  Via \M, ${\tilde C}$ has degree 19 and arithmetic genus 12. It
follows easily that $C={\tilde C}$, i.\hskip0em e.\ $C$ is cut out
scheme-theoretically
by quadrics.

Next, to see that $C$ is \aCM, note that $C$ is non-special since the
projective
dimension of the embedding system is 7, is linearly normal by construction,
and is quadratically normal by Riemann-Roch and $h^0(I_C(2))=9$ found by \M.
This
suffices to show that $C$ is \aCM\ by \cite[P.\ 222]{ACGH} or the argument in
the proof of Theorem 1.2.7 in \cite{L}.

\bigskip\noindent
{\em Proof of Proposition~\ref{main}:\/} The key point is to show that the
conditions ``\aCM'' and ``scheme-theoretically cut out by quadrics'' are dense.

A curve is \aCM\ if and only if it is projectively normal.  So $C$ is \aCM\ if
and
only if $H^1(I_C(n))=0$ for all $n\ge 0$, where $I_C$ is the ideal sheaf of
$C$.
By \cite{GLP}, $H^1(I_C(n))=0$ for all $n\ge 13$, so there are only finitely
many cohomology groups that are required to vanish in addition.  By upper
semicontinuity of $h^1(I_C(n))=\dim\ H^1(I_C(n))$, this is a Zariski open
condition in the Hilbert scheme.

As to the condition of being scheme-theoretically cut out by quadrics, we may
restrict to considering curves which are \aCM.  Let $V$ be the 9 dimensional
space
of quadrics containing $C$.  Consider the maps

\begin{equation}
\label{Vspan}
V\otimes H^0(\O_{\P7}(k))\to H^0(I_C(k+2))
\end{equation}

$V$ cuts out $C$ scheme-theoretically if and only if (\ref{Vspan}) is
surjective for some $k\ge 12$ (since $C$ is 14-regular by \cite{GLP}; a
smaller bound for effective $k$ can be given if desired).  This is again an
open condition.

Finally, let ${\rm Hilb}_{19n-11}^0$ be the subset of the Hilbert scheme
parametrizing smooth, irreducible curves in \P7 of degree 19 and genus 12. It
is open in the Hilbert scheme by \cite[P.\ 99]{GIT}.  ${\rm Hilb}_{19n-11}^0$
is
defined over Spec ${\bf Z}$ and is irreducible (its geometric fibers are
equidimensional and irreducible; this follows from the
irreducibility of ${\cal M}_{12}$ in arbitrary characteristic  \cite{DM}, and
the
non-speciality of $|L|$).

Hence the set of smooth \aCM\ curves scheme-theoretically cut out by quadrics
is
non-empty and open, hence dense, in ${\rm Hilb}_{19n-11}^0$.  This completes
the
proof of Proposition~\ref{main}.  \QED

\bigskip

It seems appropriate to conclude with some related questions.

In \cite{ACGH}, \cite[\S 3]{GL}, it was proven that a general linear system of
degree
$d\ge [(3g+4)/2]$ on a
curve $C$ of genus $g$ embeds $C$ as an \aCM\ curve.  Rather than looking for
a bound for {\em all} curves, instead one can ask:

\bigskip\noindent
Problem: Find the smallest possible $d(g)$ such that for all $d\ge d(g)$,
a general curve of genus $g$
admits a degree $d$ complete embedding which is \aCM.

\bigskip\noindent
{\em Remark.}  Suppose that $d\ge (2g+1+\sqrt{8g+1})/2$.  Then the
general degree $d$ embedding of a general curve of genus $g$ is \aCM\
\cite{BE}.
This bound is in fact sharp for {\em non-special} embeddings.
The inequality is just the solution of the inequality
$h^0(\O_{\P r}(2))\ge h^0(\O_C(2))$ for a general non-special embedding.

\bigskip
Similarly, one can ask

\bigskip\noindent
Problem: Find the smallest possible $d'(g)$ such that a general degree $d$
embedding of a general curve of genus $g$ is scheme theoretically cut out by
quadrics if $d\ge d'(g)$.

\bigskip
By work of Green and Lazarsfeld \cite[Prop.\ 2.4.2]{L}, \mbox{$d'(g)\le
[(3g+6)/2]$}, and Proposition~\ref{main} shows that this is not sharp.

\bigskip\noindent
Question: Is Proposition~\ref{main} true without restriction on the
characteristic?  Is there a counterexample to the main question with $r=6$?

\end{document}